\documentclass[11pt]{article}
\usepackage[matrix,arrow,curve,cmtip]{xy}
\usepackage{graphicx}
\usepackage{amssymb}
\usepackage{amsmath}
\usepackage{latexsym}
\usepackage{eucal}
\usepackage{theorem}
\usepackage[a4paper,width=170mm,height=230mm]{geometry}
\usepackage{enumitem}
\usepackage{tikz-cd}
\usepackage{titlesec}
\titleformat{\section}[hang]%
{\bf\filcenter\large}{\thesection.}{1ex}{}%
\titleformat{\subsection}[runin]%
{\bfseries\normalsize}{}{0ex}{}
\titleformat{\paragraph}[runin]%
{\bfseries\normalsize}{}{0ex}{}[.]

\def\to{\mbox{$\xymatrix@1@C=5mm{\ar@{->}[r]&}$}}
\def\tto{\mbox{$\xymatrix@1@C=5mm{\ar@{=>}[r]&}$}}
\def\biar{\mbox{$\xymatrix@1@C=5mm{\ar@<1.5mm>[r]\ar@<-0.5mm>[r]&}$}}
\def\adjar{\mbox{$\xymatrix@1@C=5mm{\ar@<1.5mm>@{<-}[r]\ar@<-0.5mm>[r]&}$}}
\def\iso{\mbox{$\xymatrix@1@C=6mm{\ar@{->}[r]^{\sim}&}$}}
\def\distsign{\begin{picture}(0,0)\put(0,0){\circle{4}}\end{picture}}
\def\dist{\mbox{$\xymatrix@1@C=5mm{\ar@{->}[r]|{\distsign}&}$}}

\newtheorem{theorem}{Theorem}[section]
\newtheorem{lemma}[theorem]{Lemma}

\newtheorem{corollary}[theorem]{Corollary}
{\theorembodyfont{\upshape}\newtheorem{example}[theorem]{Example}}
{\theorembodyfont{\upshape}}
\newcommand{\proof}{\noindent {\it Proof\ }: }

\def\endofproof{$\mbox{ }\hfill\Box$\par\vspace{1.8mm}\noindent}

\def\Pos{\mathsf{Pos}}
\def\Met{\mathsf{Met}}

\def\down{\ \downarrow\!}
\def\obj{\mathrm{obj}}
\def\:{\colon}
\def\1{\mathbf{1}}
\def\2{\mathbf{2}}
\def\Dist{{\sf Dist}}
\def\Cat{{\sf Cat}}
\def\Sup{{\sf Sup}}
\def\Qtld{{\sf Qtld}}
\def\Q{\mathcal{Q}}
\def\C{\mathcal{C}}
\def\D{\mathcal{D}}
\def\V{\mathcal{V}}
\def\P{\mathcal{P}}
\def\bbA{\mathbb{A}}
\def\bbB{\mathbb{B}}

\def\bbP{\mathbb{P}}
\def\bbT{\mathbb{T}}
\def\eqref#1{(\ref{#1})}
\def\inv{^{-1}}
\def\slice#1{%
\setlength{\unitlength}{1ex}
\begin{picture}(2.5,1)(0,0)
\put(-0.2,-0.5){\mbox{$/$}}
\put(1,-1){\mbox{${#1}$}}
\end{picture}}

\title{When is $\Cat(\Q)$ cartesian closed?}
\author{Isar Stubbe\footnote{Universit\'e du Littoral-C\^ote d'Opale, Calais, France; {\tt isar.stubbe@univ-littoral.fr}, corresponding author.} \quad and \quad Junche Yu\footnote{School of Mathematical Sciences, Ocean University of China, Qingdao, China; {\tt juncheyu@ouc.edu.cn}}}
\date{28 December 2025\footnote{First version posted on 7 January 2025; this revision takes referee's comments into account.}}

\begin{document}

\maketitle

\begin{quote}{\small
{\bf Abstract.} We give an elementary characterization of those quantaloids $\Q$ for which the category $\Cat(\Q)$ of $\Q$-enriched categories and functors is cartesian closed. We then unify several known cases (previously proven using \textit{ad hoc} methods) and we give some new examples. 
}\end{quote}

\section{Introduction}\label{A}

An object $A$ in a category $\C$ with finite products is said to be \textit{exponentiable} when the functor $-\times A\:\C\to \C$ has a right adjoint (whose action on an object $B$ in $\C$ is then written as $B^A$). The category $\C$ is \textit{cartesian closed} if it has all finite products and all of its objects are exponentiable. This is an important condition, e.g.\ it is a crucial axiom for toposes. But also non-topos categories can be cartesian closed: a well-known example is the category $\Pos$ of ordered sets and order-preserving maps. On the other hand, the familiar category $\Met$ of (generalized) metric spaces and non-expanding maps has all products but is not cartesian closed \cite{CH06}. These two examples are instances of the same general construction, namely categories enriched in a quantale \cite{L73}: for ordered sets the base quantale is the Boolean algebra $(\{0,1\},\wedge,\top)$, whereas for metric spaces it is Lawvere's quantale of non-negative real numbers $([0,\infty],+,0)$. In fact, quantales themselves are precisely quantaloids with a single object (like monoidal categories vs.\ bicategories), and quantaloidal enrichment -- even though slighly more involved -- captures many more useful examples, e.g.\ partial metric spaces. Thus the question arises: which property of a base quantaloid $\Q$ makes the category of $\Q$-enriched categories (which always has products) cartesian closed? In this paper, we shall give a full answer to this question. 

To be more precise, recall that a quantaloid $\Q$ is a 2-category whose homs are suplattices and in which composition distributes over suprema:
\[g\circ\big(\bigvee_if_i\Big)=\bigvee_i(g\circ f_i)\quad\mbox{ and }\quad\Big(\bigvee_ig_i\Big)\circ f=\bigvee_i(g_i\circ f).\]
Henceforth we shall always assume that $\Q$ is small, and we write $\Q_0$ for its set of objects and $\Q_1$ for its set of morphisms. A $\Q$-enriched category $\bbA$ is a set $\bbA_0$ of `objects' together with a `type' function $\bbA_0\to\Q_0\:a\mapsto ta$ and a `hom' function $\bbA_0\times\bbA_0\to\Q_1\:(a',a)\mapsto\bbA(a',a)$ such that, for all $a,a',a''\in\bbA_0$,
\[\bbA(a',a)\:ta\to ta'\mbox{ in }\Q\quad\mbox{ , }\quad \bbA(a'',a')\circ\bbA(a',a)\leq\bbA(a'',a)\quad\mbox{ and }\quad 1_{ta}\leq\bbA(a,a).\]
A $\Q$-enriched functor $F\:\bbA\to\bbB$ is an `object function' $\bbA_0\to\bbB_0\:a\mapsto Fa$ such that, for all $a,a'\in\bbA_0$,
\[t(Fa)=ta\quad\mbox{ and }\quad  \bbA(a',a)\leq\bbB(Fa',Fa).\]
With obvious composition and identities, $\Q$-categories and functors form a (large) category $\Cat(\Q)$. The terminal $\Q$-category $\bbT$ has $\bbT_0=\Q_0$ as object set, the type function is the identity, and the hom-arrow $\bbT(Y,X)$ is the top element of $\Q(X,Y)$; it is easily seen that the type function of a $\Q$-category $\bbA$ underlies the unique functor from $\bbA$ to $\bbT$. The binary product of two $\Q$-categories $\bbA$ and $\bbB$ has object set $(\bbA\times\bbB)_0=\{(a,b)\in\bbA_0\times\bbB_0\mid ta=tb\}$ with types $t(a,b)=ta=tb$ and hom-arrows $(\bbA\times\bbB)((a',b'),(a,b))=\bbA(a',a)\wedge\bbB(b',b)$; the projections are obvious and the universal property is easy to verify.

The starting point for our investigation is the following result from the article \cite{CHS09}:
\begin{theorem}[\cite{CHS09}]\label{1}
A functor $F\:\bbA\to\bbB$ between $\Q$-enriched categories is exponentiable in $\Cat(\Q)$, i.e.\ the functor
$$-\times F\:\Cat(\Q)\slice{\bbB}\to\Cat(\Q)\slice{\bbB}$$
admits a right adjoint, if and only if the following two conditions hold:
\begin{enumerate}
\item\label{e1} for every $a,a'\in\bbA$ and $\bigvee_if_i\leq\bbB(Fa',Fa)$, $$\Big(\bigvee_if_i\Big)\wedge\bbA(a',a)=\bigvee_i\Big(f_i\wedge\bbA(a',a)\Big),$$
\item\label{e2} for every $a,a''\in\bbA$, $b'\in\bbB$, $f\leq\bbB(b',Fa)$ and $g\leq\bbB(Fa'',b')$,
$$(g\circ f)\wedge\bbA(a'',a)=\bigvee_{a'\in F\inv b'}\Big((g\wedge\bbA(a'',a'))\circ(f\wedge\bbA(a',a))\Big).$$
\end{enumerate}
\end{theorem}
Since a $\Q$-category $\bbA$ is an exponentiable object of $\Cat(\Q)$ precisely when the unique functor from $\bbA$ into the terminal $\Q$-category is an exponentiable morphism, the above implies:
\begin{corollary}[\cite{CHS09}]\label{B}
A $\Q$-category $\bbA$ is an exponentiable object of $\Cat(\Q)$ if and only if
\begin{enumerate}
\item for all $a,a'\in\bbA$, $\bbA(a',a)$ is exponentiable in the suplattice $\Q(ta,ta')$,
\item for all $a,a''\in\bbA$ and $f\:ta\to Y$, $g\:Y\to ta''$ in $\Q$, 
$$(g\circ f)\wedge\bbA(a'',a)=\bigvee\left\{(g\wedge\bbA(a'',a'))\circ(f\wedge\bbA(a',a))\mid a'\in\bbA,\ ta'=Y\right\}.$$
\end{enumerate}
\end{corollary}
In particular, any locale $L$ (i.e.\ a complete lattice in which finite infima distribute over arbitrary suprema, or equivalently, a suplattice in which each element is exponentiable) can be viewed as a quantale with binary infimum as composition---and therefore as a one-object quantaloid $\Q_L$. With the Corollary above, it is easily verified that \emph{every} $\Q_L$-category is exponentiable; that is, $\Cat(\Q_L)$ is cartesian closed. 

Furthermore, the following particular situation was also pointed out in \cite{CHS09} (rephrased slightly here):
\begin{example}[\cite{CHS09}]\label{C}
Suppose that the base quantaloid $\Q$ is locally localic (meaning that each hom-suplattice is a locale) and satisfies the interchange law: for any arrows
$$\xymatrix@=12mm{
A\ar@<1mm>[r]^f\ar@<-1mm>[r]_h & B\ar@<1mm>[r]^g\ar@<-1mm>[r]_k & C}$$
in $\Q$, the equality $(g\circ f)\wedge(k\circ h)=(g\wedge k)\circ(f\wedge h)$ holds. Then a $\Q$-category $\bbA$ is an exponentiable object of $\Cat(\Q)$ if and only if, for all $a,a''\in\bbA$ and $Y\in\Q$,
\begin{equation*}\label{D}
(\top_{Y,ta''}\circ\top_{ta,Y})\wedge\bbA(a'',a)=\bigvee\left\{\bbA(a'',a')\circ\bbA(a',a)\mid a'\in\bbA,\ ta'=Y\right\}.
\end{equation*}
\end{example}
However, \cite[Example 5.2]{CHS09} then contains the erroneous claim that all {\em free quantaloids}\footnote{Given a (small) category $\C$, the free quantaloid $\P\C$ has the same objects as $\C$, the hom-suplattice $\P\C(X,Y)$ is the powerset $\P(\C(X,Y))$, composition in $\P\C$ is done ``elementwise'',
$$\mbox{for }S\subseteq\C(X,Y)\mbox{ and }T\subseteq\C(Y,Z)\mbox{ define }T\circ S=\{t\circ s\mid s\in S,t\in T\},$$
and the identity on an object $X$ is $\{1_X\}$. If $\C$ has only one object -- so it is effectively a monoid -- then this construction produces the \emph{free quantale} on that monoid. See \cite{R96}.} satisfy the interchange law; and from this it is then deduced that the category of categories enriched in a \textit{free quantale} is always cartesian closed. The mistake is easily recognized: a quantale $Q=(Q,\bigvee,\circ,1)$ satisfying the interchange law has a second monoid structure $(Q,\wedge,\top)$ which, by the ``Eckmann-Hilton argument'', is necessarily identical to $(Q,\circ,1)$ -- in other words, such a quantale $Q$ is nothing but a locale (it has $\circ=\wedge$) -- yet this is not the case for a non-trivial free quantale! So this leaves it an open question whether or not $\Cat(\Q)$ is cartesian closed whenever $\Q$ is a free quantaloid---and this instigated our research presented in this paper.

In the following section, a series of Lemmas builds up to our main Theorem \ref{QCatCC} that gives an elementary necessary-and-sufficient condition on a quantaloid $\Q$ for the category $\Cat(\Q)$ to be cartesian closed. With this characterization, we then unify several known cases (previously proven using \textit{ad hoc} methods), and we give some new examples. Our final Example \ref{FreeQ} corrects the erroneous claim from \cite{CHS09}: we point out that $\Cat(Q)$ is \emph{never} cartesian closed for a non-trivial free quantale $Q$.

\section{Cartesian closedness of $\Cat(\Q)$}

Corollary \ref{B} shows exactly how the exponentiability of each individual $\Q$-category $\bbA$ depends on the base quantaloid $\Q$. Thus, to find (at least) \emph{necessary} conditions on $\Q$ for $\Cat(\Q)$ to be cartesian closed, we can put different ``test'' categories in place of $\bbA$. A particularly simple kind of $\Q$-category will be useful: for any arrow $f\:X\to Y$ in $\Q$, let $\bbP_f$ denote the $\Q$-category with two (different) objects, $*_1$ and $*_2$, with types $t*_1=X$ and $t*_2=Y$, and whose homs are
\[\bbP_f(*_1,*_1)=1_X,\quad\bbP_f(*_2,*_2)=1_Y,\quad\bbP_f(*_2,*_1)=f,\quad\bbP_f(*_1,*_2)=\perp_{Y,X}.\]
This category $\bbP_f$ is the \emph{collage} of the one-element distributor $(f)\:\1_X\dist\1_Y$, which itself is the image of the arrow $f\:X\to Y$ under the inclusion of quantaloids $\Q\to\Dist(\Q)$ (the codomain of which is the quantaloid of $\Q$-categories and $\Q$-distributors).
\begin{lemma}\label{CCimplieslocallylocalic}
If $\Cat(\Q)$ is cartesian closed, then $\Q$ is locally localic.
\end{lemma}
\proof
For any arrow $f\:X\to Y$ in $\Q$, by hypothesis we know that $\bbP_f$ is exponentiable in $\Cat(\Q)$, which by the first condition in Corollary \ref{B} implies exponentiability of $f$ in $\Q(X,Y)$. Thus each hom-suplattice of $\Q$ must indeed be a locale.
\endofproof
\begin{lemma}\label{composition-bottom}
If $\Cat(\Q)$ is cartesian closed, then for any 
$\begin{tikzcd}[baseline=(W.base),every label/.append style = {font = \normalsize},row sep=0ex]  
X \ar[r,"f"]&Y\ar[r,"g"]& |[alias=W]| Z
\end{tikzcd}$ 
with $X\neq Y\neq Z$ in $\Q$ we have that $g\circ f=\bot_{X,Z}$.
\end{lemma}
\proof 
For any objects $X$ and $Z$ in $\Q$, consider the top arrow $\top_{X,Z}\:X\to Z$; we know by hypothesis that the $\Q$-category $\bbP_{\top_{X,Z}}$ is exponentiable in $\Cat(\Q)$. Given any $f\colon X\to Y$ and $g\colon Y\to Z$ in $\Q$ with $X\neq Y\neq Z$,
the absence of objects of type $Y$ in $\bbP_{\top_{X,Z}}$ makes the right hand side of the second condition in Corollary \ref{B} the supremum of the empty subset of $\Q(X,Z)$. Thus 
\[g\circ f = (g\circ f)\land\top_{X,Z}=(g\circ f)\land\bbP_{\top_{X,Z}}(*_2,*_1)=\bigvee\emptyset=\bot_{X,Z}\]
as claimed.
\endofproof
\begin{lemma}\label{X=YneqZ}
If $\Cat(\Q)$ is cartesian closed, then
\begin{enumerate}[label=\rm(\roman*)]
\item\label{cond1} for any 
$\begin{tikzcd}[baseline=(W.base),every label/.append style = {font = \normalsize},row sep=0ex]  
X \ar[r,"f"]&X\ar[r,"g",shift left=1.5]\ar[r,"h" below,shift right] & |[alias=W]| Y
\end{tikzcd}$ 
with $X\neq Y$ in $\Q$ we have that $(g\circ f)\wedge h = (g\land h)\circ(f\land 1_X)$;
\item\label{cond2} for any 
$\begin{tikzcd}[baseline=(W.base),every label/.append style = {font = \normalsize},row sep=0ex]  
X \ar[r,"f",shift left=1.5]\ar[r,"h" below,shift right]&Y\ar[r,"g"] & |[alias=W]| Y
\end{tikzcd}$ 
with $X\neq Y$ in $\Q$ we have that $(g\circ f)\wedge h = (g\land 1_Y)\circ(f\land h)$.
\end{enumerate}
\end{lemma}
\proof
We prove the first of these conditions; the second is similar. By hypothesis we know that $\bbP_h$ is exponentiable in $\Cat(\Q)$, so we can compute with the second condition in Corollary \ref{B} that 
\[(g\circ f)\wedge h 
	= (g\circ f)\wedge\bbP_h(*_2,*_1)
	= (g\wedge\bbP_h(*_2,*_1))\circ(f\wedge\bbP_h(*_1,*_1))
	= (g\wedge h)\circ(f\wedge 1_X),\]
as wanted.
\endofproof
\begin{lemma}\label{X=Y=Z}
If $\Cat(\Q)$ is cartesian closed, then for any 
$\begin{tikzcd}[baseline=(W.base),every label/.append style = {font = \normalsize},row sep=0ex]  
X \ar[r,"f",shift left=0.5,bend left]\ar[r,"h" below,shift right=0.5,bend right]\ar[r,"g" description] & |[alias=W]| X
\end{tikzcd}$
in $\Q$ we have that 
\[(g\circ f)\land h = ((g\land h)\circ (f\land 1_X))\lor ((g\land 1_X)\circ (f\land h)).\]
\end{lemma}
\proof
As before, we use the hypothetical exponentiability of $\bbP_h$ to compute with the second condition in Corollary \ref{B} that
	\begin{align*}
	(g\circ f)\wedge h 
	& = (g\circ f)\wedge\bbP_h(*_2,*_1)\\ 
	& = ((g\land \bbP_h(*_2,*_1))\circ (f\land \bbP_h(*_1,*_1)))\lor ((g\land \bbP_h(*_2,*_2))\circ (f\land \bbP_h(*_2,*_1)))\\
	& = ((g\land h)\circ (f\land 1_X))\lor ((g\land 1_X)\circ (f\land h)),
	\end{align*}
noting that now $\bbP_h$ has two objects of type $X$.
\endofproof
As it now turns out, the necessary conditions for $\Cat(\Q)$'s cartesian closedness established in the four previous Lemmas, are also sufficient:
\begin{theorem}\label{QCatCC}
Let $\Q$ be a small quantaloid. The category $\Cat(\Q)$ is cartesian closed if and only if
\begin{enumerate}[label=\rm(\roman*)]
\item\label{ll} $\Q$ is locally localic, and
\item\label{thecondition} for any 
$\begin{tikzcd}[baseline=(W.base),every label/.append style = {font = \normalsize},row sep=0ex,column sep=3ex]
   & Y \ar[dr,bend left, "g"] \\   
  X \ar[rr,"h" below]\ar[ur,bend left,"f"] && |[alias=W]| Z
\end{tikzcd}$ 
in $\Q$ we have that
\[(g\circ f)\land h =
\begin{cases}
\bot_{X,Z} & \mbox{ if }X\neq Y\neq Z,\\
(g\land h)\circ (f\land 1_X) & \mbox{ if }X = Y\neq Z,\\
(g\land 1_Z)\circ (f\land h) & \mbox{ if }X\neq Y = Z,\\
\Big((g\land h)\circ (f\land 1_X)\Big)\lor \Big((g\land 1_Z)\circ (f\land h)\Big) & \mbox{ if }X=Y=Z.
\end{cases}\]
\end{enumerate}
\end{theorem}
\proof
The necessity of the conditions in this statement follows from Lemmas \ref{CCimplieslocallylocalic}, \ref{composition-bottom}, \ref{X=YneqZ} and \ref{X=Y=Z}. For the sufficiency, suppose that $\Q$ satisfies both conditions in the statement. Now let $\bbA$ be any $\Q$-category; we shall verify that both conditions in Corollary \ref{B} hold. Since $\Q$ is locally localic, the first condition in Corollary \ref{B} is certainly satisfied. As for the second condition in Corollary \ref{B}, we only need to verify that the left hand side is less than or equal to the right hand side, because the converse inequality is always true (using the composition law in $\bbA$). Let $a,a''\in\bbA$ and $f\colon ta\to Y$, $g\colon Y\to ta''$ in $\Q$. We distinguish four cases:
\begin{itemize}
	\item If $ta\neq Y\neq ta''$, then
	\[(g\circ f)\wedge\bbA(a'',a)=\bot_{ta,ta''}=\bigvee\left\{(g\wedge\bbA(a'',a'))\circ(f\wedge\bbA(a',a))\mid a'\in\bbA,\ ta'=Y\right\}\]
	because any arrow in $\Q$ that factors through an object which is neither its domain nor its codomain is the bottom arrow. (This follows from the first of the four cases in condition \ref{thecondition} by letting $h=\top_{X,Z}$.) 
	\item If $ta=Y\neq ta''$, then 
	\begin{align*}
	(g\circ f)\wedge\bbA(a'',a)
	& = (g\land\bbA(a'',a))\circ (f\land 1_{ta})\\
	& \leq (g\land \bbA(a'',a))\circ (f\land \bbA(a,a))\\
	& \leq \bigvee\left\{(g\wedge\bbA(a'',a'))\circ(f\wedge\bbA(a',a))\mid a'\in\bbA,\ ta'=Y\right\}.
	\end{align*}
	where we use the second case of condition \ref{thecondition} and the fact that $1_{ta}\leq\bbA(a,a)$.
	\item If $ta\neq Y = ta''$, then 
	\begin{align*}
	(g\circ f)\wedge\bbA(a'',a) 
	& = (g\land 1_{ta''})\circ (f\land \bbA(a'',a))\\
	& \leq (g\land \bbA(a'',a''))\circ (f\land \bbA(a'',a))\\
	& \leq \bigvee\left\{(g\wedge\bbA(a'',a'))\circ(f\wedge\bbA(a',a))\mid a'\in\bbA,\ ta'=Y\right\}
	\end{align*}
		where we use the third case of condition \ref{thecondition} and the fact that $1_{ta''}\leq\bbA(a'',a'')$.
	\item If $ta=Y=ta''$, then
	\begin{align*}
	(g\circ f)\wedge\bbA(a'',a) 
	& = ((g\land \bbA(a'',a))\circ (f\land 1_{ta}))\lor ((g\land 1_{ta''})\circ (f\land \bbA(a'',a)))\\
	& \leq ((g\land \bbA(a'',a))\circ (f\land \bbA(a,a)))\lor ((g\land \bbA(a'',a''))\circ (f\land \bbA(a'',a)))\\
	& \leq \bigvee\left\{(g\wedge\bbA(a'',a'))\circ(f\wedge\bbA(a',a))\mid a'\in\bbA,\ ta'=Y\right\}
	\end{align*}
	now using the fourth case of condition \ref{thecondition} and both $1_{ta}\leq\bbA(a,a)$ and $1_{ta''}\leq\bbA(a'',a'')$.
	\end{itemize}
This concludes the proof.
\endofproof
To end this section, let us remark that, if $\Cat(\Q)$ is cartesian closed, then surely, for every $f\:X\to Y$ in $\Q$, the $\Q$-category $\bbP_f$ is exponentiable. However, the converse is true too: the proofs of Lemmas \ref{CCimplieslocallylocalic}, \ref{composition-bottom}, \ref{X=YneqZ} and \ref{X=Y=Z} only use the exponentiability of categories of the form $\bbP_f$, yet those Lemmas establish the (necessary and) sufficient conditions in Theorem \ref{QCatCC} for $\Cat(\Q)$ to be cartesian closed.

\section{Examples}

\begin{example}[Quantales, locales]\label{quantales}
Let $Q=(Q,\circ,1)$ be a quantale. Then $\Cat(Q)$ is Cartesian closed if and only if the underlying suplattice of $Q$ is a locale and 
\begin{equation}\label{cond-quantale}
\mbox{ for all $a,b,c\in Q$:}\quad(a\circ b)\wedge c = ((a\wedge c)\circ (b\wedge 1))\vee((1\wedge a)\circ(b\wedge c)),
\end{equation}
because the other conditions in Theorem \ref{QCatCC} are void. Clearly, this condition is met whenever the multiplication in $Q$ is in fact the binary infimum, i.e.\ when $Q$ is nothing but a locale.
\end{example}

\begin{example}[A non-locale example]\label{nonlocale}
Endow the set $Q=\{0,\frac{1}{2},1\}$ with the natural order and the multiplication $x\circ y=\max\{x+y-1,0\}$; this is exactly the truth-value table of (the conjunction in) {\L}ukasiewicz's three-valued logic \cite{L20}. It is easy to check that $Q$ is a quantale that satisfies the conditions in Theorem \ref{QCatCC}, hence $\Cat(Q)$ is cartesian closed. This shows that $\Cat(Q)$ may be cartesian closed even when $Q$ is not a locale, that is, it may have $\circ\neq\land$ (yet its underlying suplattice is necessarily a locale). This example appears in \cite[Example 4.8]{LZ16} with an \textit{ad hoc} proof.
\end{example}

\begin{example}[Integral quantales, $t$-norms]
For the quantale $Q=\{0 < 1 < \top\}$ (with the only possible multiplication that has $1$ as neutral element), $\Cat(Q)$ is cartesian closed: this shows that condition \eqref{cond-quantale} does not imply that $1=\top$ in $Q$. However, if $Q$ is an integral quantale (meaning that $1=\top$), then the condition in \eqref{cond-quantale} simplifies to
\begin{equation}\label{cond-int-quantale}
\mbox{ for all $a,b,c\in Q$:}\quad(a\circ b)\wedge c = ((a\wedge c)\circ b)\vee(a\circ(b\wedge c)).
\end{equation}
Recall that a \emph{left-continuous $t$-norm} is exactly an integral quantale whose underlying suplattice is the real interval $[0,1]$ (with natural order); such a $t$-norm is \emph{continuous} whenever multiplication is a continuous function (in each variable). As the underlying suplattice of a (left-)continuous $t$-norm $Q=([0,1],\circ,1)$ is a locale, $\Cat(Q)$ is cartesian closed if and only if \eqref{cond-int-quantale} holds. Theorem 4.7 of \cite{LZ16} implies that the only \emph{continuous} $t$-norm satisfying \eqref{cond-int-quantale} is the G\"odel $t$-norm, that is, whose multiplication is given by binary infimum. However, there are other \emph{left-continuous} $t$-norms satisfying \eqref{cond-int-quantale}, for example when putting 
\[x\circ y=\left\{\begin{array}{ll}
0&\mbox{ if }x,y\leq\frac{1}{2}, \\
x\wedge y&\mbox{ otherwise.}
\end{array}\right.\]
As a final remark, condition \eqref{cond-int-quantale} was first observed in \cite[Theorem 4.6]{LL25}, but only for so-called complete subquantales of continuous $t$-norms (such as Example \ref{nonlocale} above), and with a different proof than ours.
\end{example}

\begin{example}[Squares, cubes and idempotents]
If $Q$ is an integral quantale satisfying \eqref{cond-int-quantale}, then an easy computation shows that (writing $aa$ for $a\circ a$)
\[aa=aa\wedge aa=(a\wedge aa)a\vee a(a\wedge aa)=aaa\]
for any $a\in Q$. (This, by the way, clearly does not hold in the Lawvere quantale $([0,\infty],\bigwedge,+,0)$, so the category of generalized metric spaces is not cartesian closed, cf.\ \cite{CH06}.) It follows that every square in such a $Q$ is idempotent, and in \cite[Theorem 4.6]{LL25} it is shown that, for localic subquantales of continuous $t$-norms, this fact is equivalent to \eqref{cond-int-quantale}.
\end{example}

\begin{example}[A quantaloidal example]
Let $\Q$ be the split-idempotent completion of the two-element Boolean algebra $\2$, that is, $\Q$ has objects and arrows as in
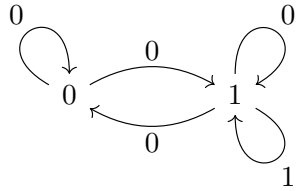
\[\begin{tikzcd}[column sep=10ex,every label/.append style = {font = \normalsize}]
0\ar[r,bend left,"0"]\ar[out=150, in=90, loop, distance=3em,"0"] & 1\ar[l,bend left,"0"]\ar[out=90, in=30, loop, distance=3em,"0"]\ar[out=-30, in=-90, loop, distance=3em,"1"]
\end{tikzcd}\]
with composition given by binary infimum. This quantaloid satisfies the conditions in Theorem \ref{QCatCC}, so $\Cat(\Q)$ is cartesian closed. Remark that $\Q$ is precisely the quantaloid of diagonals in $\2$, \textit{viz.}\ $\Q=\D(\2)$, so that $\Cat(\Q)$ is exactly the category of \emph{partial $\2$-enriched categories} \cite{HS18}.
\end{example}

\begin{example}[Coproducts of quantaloids]
Given two small quantaloids $\Q_1$ and $\Q_2$, their coproduct in the category $\Qtld$ of quantaloids and homomorphisms (i.e.\ the category of $\Sup$-enriched categories and $\Sup$-enriched functors), that we shall denote by $\Q$, is constructed as follows:
\begin{itemize}
\item $\obj(\Q)=\obj(\Q_1)\uplus\obj(\Q_2)$,
\item $\Q(X,Y)=\left\{\begin{array}{ll}
\Q_1(X,Y)&\mbox{ if }X,Y\in\obj(Q_1),\\
\Q_2(X,Y)&\mbox{ if }X,Y\in\obj(Q_2),\\
\{\bot\} & \mbox{ otherwise},
\end{array}\right.$
\end{itemize}
and composition and identities are defined in the obvious way. It is straightforward to verify with Theorem \ref{QCatCC} that both $\Cat(\Q_1)$ and $\Cat(\Q_2)$ are cartesian closed if and only if $\Cat(\Q)$ is cartesian closed. This generalizes the previous example, which can be interpreted as the coproduct of two (one-object suspensions of) locales.
\end{example}

\begin{example}[Another quantaloidal example]
Let $L$ be a locale, pick two elements $u,v\in L$, and consider the quantaloid $\Q$ with
\begin{itemize}
\item $\obj(\Q)=\{u,v\}$,
\item $\Q(u,u)=\down u$, $\Q(v,v)=\down v$, $\Q(u,v)=\down(u\wedge v)$ and $\Q(v,u)=\bot$,
\end{itemize}
with composition given by binary infima. This quantaloid satisfies the conditions in Theorem \ref{QCatCC}, so $\Cat(\Q)$ is cartesian closed. If $u\wedge v=\bot$, then this $\Q$ is the coproduct in $\Qtld$ of $\down u$ and $\down v$ seen as one-object quantaloids; otherwise it is not a coproduct of one-object quantaloids.
\end{example}

\begin{example}[Diagonals]
Let $L$ be a locale; we already saw that $\Cat(L)$ is cartesian closed. Now let $\D(L)$ be the quantaloid of diagonals in $L$ (equivalently in this case, obtained by splitting idempotents in $L$). It is locally localic but as soon as $L$ has at least three elements, $\Cat(\D(L))$ is not cartesian closed because the requirement in Lemma \ref{composition-bottom} does not hold. This shows that the properties in Theorem \ref{QCatCC} are not stable under the splitting of idempotents nor under the construction of diagonals (see \cite[Example 2.14]{S14} for more on diagonals in a quantaloid). Incidentally, both $L$ and $\D(L)$ satisfy the interchange law, so this also goes to show that this is neither a necessary nor a sufficient condition for cartesian closedness.
\end{example}

\begin{example}[Free quantaloids]\label{FreeQ}
The free quantaloid $\P\C$ on a small category $\C$ is always locally localic: its local suprema/infima are unions/intersections in powersets. It furthermore satisfies condition \ref{thecondition} in Theorem \ref{QCatCC} (i.e.\ $\Cat(\P\C)$ is cartesian closed) if and only if
\begin{equation}\label{cond-free-qtld}
\mbox{in the category $\C$, when two morphisms compose then at least one of them is an identity.}
\end{equation}
(Put differently, this is a category in which \textit{every morphism is prime}.) Indeed, \eqref{cond-free-qtld} says in particular that the only endomorphisms in $\C$ are identities, which makes it straightforward to check the four cases in Theorem \ref{QCatCC}--\ref{thecondition}. Conversely, first consider an endomorphism $f\:X\to X$ in $\C$, and suppose that $f\neq 1_X$. Putting $F=G=\{f\}$ and $H=\C(X,X)$ in $\P\C(X,X)$, we get from the fourth case in Theorem \ref{QCatCC}--\ref{thecondition} that 
\[\{f\circ f\}=(G\circ F)\cap H=\Big((G\cap H)\circ (F\cap\{1_X\})\Big)\cup\Big(G\cap\{1_X\})\circ(F\cap H)\Big)=\emptyset\cup\emptyset,\]
a contradiction; so all endomorphisms in $\C$ are identities. Now consider a composable pair of morphisms in $\C$, say $f\:X\to Y$ and $g\:Y\to Z$, neither of which is an identity; because all endomorphisms are identities, this implies that $X\neq Y\neq Z$. Puting $F=\{f\}\in\P\C(X,Y)$, $G=\{g\}\in\P\C(Y,Z)$ and $H=\C(X,Z)\in\P\C(X,Z)$ we find from the first case in Theorem \ref{QCatCC}--\ref{thecondition} that 
\[\{g\circ f\}=G\circ F\cap H=\emptyset,\]
a contradiction. Thus such a composable pair cannot exist in the first place.

For a one-object category, i.e.\ a monoid $M$, this shows that $\Cat(\P M)$ is cartesian closed if and only if $M=\{1\}$.
\end{example}

\section*{Acknowledgment}

The work reported in this article was carried out while Junche Yu was a postdoctoral researcher at the Universit\'e du Littoral-C\^ote d'Opale in Calais, France.

\section*{Funding statement}

No additional funding was received for this research.

\end{document}